\documentstyle[amssymb]{article}

\newcommand{\proof}{{\bf Proof:  }}
\newcommand{\remark}{{\bf Remark:  }}
\newcommand{\remarks}{{\bf Remarks:  }}
\newcommand{\example}{{\bf Example:  }}
\newcommand{\examples}{{\bf Examples:  }}

\newcommand{\dimv}{\underline{\dim}}

\newcommand{\hb}{\newline\hspace*{\fill}$\Box$}

\newcommand{\supp}{\mbox{supp }}

\newcommand{\ses}[3]
{\mbox{$0 \rightarrow #1 \rightarrow #2 \rightarrow #3 \rightarrow 0$}}

\newtheorem{theorem}{Theorem}[section]
\newtheorem{lemma}[theorem]{Lemma}
\newtheorem{definition}[theorem]{Definition}
\newtheorem{proposition}[theorem]{Proposition}
\newtheorem{corollary}[theorem]{Corollary}

\newtheorem{conjecture}[theorem]{Conjecture}

\begin{document}

\parindent0pt

\title{\bf The Harder-Narasimhan system in quantum groups
and cohomology of quiver moduli}
\author{Markus Reineke\\[2ex] BUGH Wuppertal, Gau\ss str. 20, D-42097
Wuppertal, Germany\\ (e-mail: reineke@math.uni-wuppertal.de)}
\date{}

\maketitle

\begin{abstract} Methods of Harder and Narasimhan from the theory of moduli of 
vector bundles are applied to moduli of quiver representations. Using the Hall 
algebra approach to quantum groups, an analog of the Harder-Narasimhan 
recursion is constructed 
inside the quantized enveloping algebra of a Kac-Moody algebra. This leads to 
a canonical orthogonal system, the HN system, in this algebra. Using a 
resolution of the recursion, an explicit formula for the HN system is given.
As an application, explicit formulas for Betti numbers of the cohomology of 
quiver moduli are derived, generalizing several results on the cohomology of 
quotients in 'linear algebra type' situations.
\end{abstract}

\section{Introduction}

The geometry of quiver representations is by now recognized as an area
with many connections to such diverse fields as representation theory of 
algebras, 
geometric invariant theory, quantum group theory and representation theory of 
(Kac-Moody) Lie algebras.\\[1ex]
The connection to methods of Geometric Invariant Theory is provided by the 
construction of moduli spaces of quiver representations of \cite{Ki}. In 
particular, one can expect analogies to the theory of moduli of vector bundles 
on curves, 
generalizing the approach advertised in \cite{Ra}.\\[1ex]
On the other hand, the connection to quantum group theory is given by the 
realization of quantized enveloping algebras of (Kac-Moody) Lie algebras of 
\cite{Ri}, \cite{Gr} via the Hall algebra approach, which can be interpreted (see \cite{Ka}) 
as a convolution algebra construction on parameter spaces of quiver 
representations.\\[1ex]
The aim of the present paper is to develop a synthesis of both methods. We 
start with a particular instance of the above mentioned analogy to vector 
bundle theory, namely the Harder-Narasimhan recursion \cite{HN}, which was 
originally used for computing Betti numbers of moduli spaces.\\[1ex]
The first main result of this paper (Proposition \ref{rec}) is a 
materialization of the Harder-Narasimhan recursion in the quantized enveloping 
algebra of a symmetric Kac-Moody Lie algebra. It leads to a canonical
orthogonal system, the HN system, in such algebras (Theorem \ref{hnsystem}), 
which is recursively computable.\\[1ex]
The HN system comprises a surprising amount of information on the moduli 
spaces of (semi-)stable quiver representations: as the second main result, 
just by evaluating at a character of the quantum group, we recover in Theorem 
\ref{betti} the Betti numbers of such moduli spaces completely (modulo the 
standard assumption (see e.g. \cite{Kir}) that semistability and stability 
coincide). The proof uses the Weil conjectures, in analogy to the original 
approach of Harder and Narasimhan in the vector bundle situation (see also 
\cite{Kir}, \cite{Goe} for similar situations).\\[1ex]
The third main result is a resolution of the Harder-Narasimhan recursion, in 
the spirit of \cite{Za}, \cite{LR} in the vector bundle case (Theorem 
\ref{ex}), together with a fast algorithm for Betti number computation 
(Corollary \ref{algo}). Whereas the 
cited works use involved explicit calculations, resp.~the Langlands lemma from 
the theory of Eisenstein series, the present proof uses only some simple 
(polygonal and simplicial) combinatorics. It should be noted that our 
materialization of the Harder-Narasimhan recursion in a non-commutative 
algebra (the quantum group) is already anticipated in (\cite{Za}, p.~457), 
where one of the key insights for resolving the recursion is a noncommutative 
approach to certain polynomial expressions.\\[1ex]
At the moment, the immediate applications of the HN system to quantum group 
theory are largely conjectural: one can expect explicit descriptions of PBW 
type bases, 
and applications to the general structure theory of Hall algebras in the 
spirit of 
\cite{SV} (see section \ref{examples}).\\[1ex]
The applications of the cohomology formulas are more direct
(and numerous), in  that they unify, generalize, and make more
explicit, several known results on  cohomology of moduli of
'linear algebra type' problems, for example, the  formulas for
Betti numbers of sequences of subspaces in (\cite{Kir}, 16),
and  for families of linear maps in \cite{D}, \cite{ES}.\\[1ex]
In another direction, note the usage of some of the present
methods in 
\cite{CBV} for a partial proof of the Kac conjectures (see \cite{Kac}). 
Finally, one can hope for analogs of the present methods to be possible in the 
setup of Hall algebras for categories of coherent sheaves on curves in the 
framework of \cite{Ka}, which would probably bring the subject back to its 
roots in vector bundle theory.\\[2ex]
The paper is organized as follows: in section 2, we recall the notions of 
(semi-) stability and of the Harder-Narasimhan filtration for categories of 
quiver representations. These are applied in section 3 to the representation 
varieties of quivers; we construct the Harder-Narasimhan stratification of 
these spaces and prove some basic geometric properties. In section 4, we first 
recall the Hall algebra construction and its relation to quantum groups. Then 
we use the results of the previous sections to obtain a quantum group version 
of the Harder-Narasimhan recursion, yielding the above mentioned 
Harder-Narasimhan system. Section 5 is devoted to the resolution of the 
recursion. In section 6, we take a closer look at the geometry of quiver 
moduli over finite fields, in order to apply the Weil conjectures, and to 
derive the formulas for Betti numbers. Finally, section 7 contains some 
conjectural applications, and the above mentioned examples of quiver moduli 
which allow for explicit cohomology formulas.\\[2ex] 
This work was done while the author participated in the TMR-network ERB
RMRX-CT97-0100 "Algebraic Lie Representations". I would like to thank K. 
Bongartz, R. Huber and M. 
Rapoport for discussions on the methods of section \ref{cqm}, E. Ossa for a 
simplification in the proof of Lemma \ref{ossa}, and B. Medeke for inspiring 
discussions on the algorithm \ref{algo}. Moreover, I would like to thank L. 
Hille, S. K\"onig, S. Orlik, A. Schofield, M. Van den Bergh and
E. Vasserot for  interesting remarks on the results of this
paper.

\section{Quiver representations and the HN filtration}\label{hnf}

Let $Q$ be a finite quiver (oriented graph) without oriented cycles. We denote 
by $I$ its
set of vertices, and by $r_{ij}$ for $i,j\in I$ the number of arrows from $i$ 
to $j$. The free abelian group ${\bf Z}I$ generated by $I$ carries
a bilinear form, the Euler form, defined by $\langle 
i,j\rangle=\delta_{ij}-r_{ij}$.\\[1ex]
Let ${\bf k}$ be an arbitrary field. We consider the well-known category 
$\bmod_{\bf k}Q$
of finite-dimensional ${\bf k}$-representations of $Q$ (see \cite{ARS} for 
general notions and facts on quiver representations). For a representation 
$X$ 
given by ${\bf k}$-vector spaces $X_i$ for $i\in I$, we define its dimension 
type $\dimv X$
as $\sum_{i\in I}\dim_{\bf k}(X_i)i\in{\bf N}I$. We have
$$\langle\dimv X,\dimv Y\rangle=\dim_{\bf k}{\rm Hom}_{{\bf k}Q}(X,Y)-\dim_{\bf k}{\rm 
Ext}^1_{{\bf k}Q}(X,Y)$$
for all representations $X,Y\in\bmod_{\bf k}Q$.\\[1ex]
We introduce a (non-canonical) notion of stability in
$\bmod_{\bf k}Q$. Fix once and for all a linear form $\Theta=\sum_{i\in
I}\Theta_ii^*$ on ${\bf Z}I$, called a  weight for $Q$. We 
define a slope function on ${\bf N}I\setminus \{0\}$ by $\mu(d)=\Theta(d)/\dim 
d$, where the linear form $\dim$ on ${\bf Z}I$ is defined by $\dim i=1$ for 
$i\in I$. Using the slope function, we can define a notion of stability for
representations of $Q$; for a representation $0\not=
X\in\bmod_{\bf k}Q$, we  set 
$\mu(X)=\mu(\dimv X)$.

\begin{definition} A ${\bf k}$-representation $X$ of $Q$ is called semistable
(resp.~stable)
if $\mu(U)\leq 
\mu(X)$ (resp.~$\mu(U)<\mu(X)$) for all proper subrepresentations 
$0\not=U\subset X$.
\end{definition}

\examples In each of the following examples, it is easy to work out the 
corresponding meaning of stability; we just state the results. The examples 
illustrate that the quiver setup unifies several 'linear algebra type' 
problems. They will be used as standard examples throughout the paper.
\begin{description}
\item[A] If $Q=i_1\rightarrow i_2\rightarrow\ldots\rightarrow i_n$ and 
$\Theta(i_k)=-k$, then the stable representations are precisely the unique 
indecomposables $E_{kl}$ of dimension type $i_k+\ldots+i_l$ for $k\leq l$, and 
the semistables are their powers (see section \ref{examples} for a possible 
generalization).  
\item[B] Let $Q$ be given by vertices $I=\{i_0,i_1,\ldots,i_n\}$ and arrows 
$i_k\rightarrow i_0$ for $k=1\ldots n$. Let $d=mi_0+i_1+\ldots+i_n$, 
$\Theta=-i_0^*$. Then a semistable (resp.~stable) representation corresponds 
to a tuple $(v_1,\ldots,v_n)$ of non-zero vectors in ${\bf k}^m$ such that for 
any 
proper non-zero subspace $U\subset {\bf k}^m$, the number of vectors $v_k$ in 
$U$ is 
$\leq$ (resp.~$<$) $\frac{n}{m}\cdot\dim U$. Thus, we arrive at one of the 
principal 
examples of Mumford's GIT (\cite{Mu}, 3.).
\item[C] More generally, fix $r,N\in{\bf N}$. Let $Q$ be given by vertices 
$I=\{i_0,i_{\nu,p}\, :\, \nu=1\ldots r,\, p=-N\ldots N\}$, arrows 
$i_{\nu,p}\rightarrow i_{\nu,p-1}$ for $\nu=1\ldots r$, $p=N\ldots 1-N$, and 
$i_{\nu,-N}\rightarrow i_0$ for $\nu=1\ldots r$. Let 
$d=\sum_{\nu,p}d_{\nu,p}i_{\nu,p}+d_0i_0$, where $0=d_{\nu,N}\leq\ldots\leq 
d_{\nu,-N}=d_0$ for $\nu=1\ldots k$, and $\Theta=-i_0^*$. Then the semistable 
(resp.~stable) representations correspond to $d_0$-dimensional ${\bf 
k}$-vector 
spaces $V$ together with a family of descending flags $F^p_\nu$ (where $\dim 
F^p_\nu=d_{\nu,p}$), fulfilling the conditions of \cite{Fa}.
\item[D] If $Q=i\stackrel{(n)}{\rightarrow}j$ (with $n$ arrows pointing from 
$i$ to $j$), $d=ai+bj$, $\Theta=i^*$, then semistable (resp.~stable) 
representations correspond to tuples $(f_1\ldots f_n)$ of linear maps 
$f_k:{\bf k}^a\rightarrow {\bf k}^b$ such that for each non-zero proper 
subspace $U\subset 
{\bf k}^a$, the dimension of $\sum_kf_k(U)$ is $\geq$ (resp.~$>$) 
$\frac{b}{a}\cdot\dim U$ (see 
\cite{ES}, \cite{D}).
\end{description}

The following properties of semistable (resp.~stable) representations are 
well-known and easy to prove (see e.g. \cite{HN}, 
\cite{Sh}, \cite{Ru}, \ldots).

\begin{lemma}\label{fl} Given a short exact sequence $\ses{M}{X}{N}$ in 
$\bmod_{\bf k}Q$, we 
have
$$\mu(M)\leq\mu(X)\mbox{ iff }\mu(X)\leq\mu(N)\mbox{ iff }\mu(M)\leq\mu(N)$$
and
$$\min(\mu(M),\mu(N))\leq\mu(X)\leq\max(\mu(M),\mu(N)).$$
If $\mu(M)=\mu(X)=\mu(N)$, then $X$ is semistable if and only if $M$ and $N$ 
are semistable.
\end{lemma}

Denote by $\bmod^{\mu}_{\bf k}Q$ the full subcategory of $\bmod_{\bf k}Q$ 
consisting of 
se\-mi\-sta\-ble representations of $Q$ of slope $\mu\in{\bf Q}$.

\begin{lemma} For all $\mu\in{\bf Q}$, the category $\bmod^{\mu}_{\bf k}Q$ is 
an 
abelian subcategory of $\bmod_{\bf k}Q$
whose simple ob\-jects are the sta\-ble representations of $Q$ of slope $\mu$. 
Moreover, we have
${\rm Hom}_Q(\bmod^{\mu}_{\bf k}Q,\bmod^{\nu}_{\bf k}Q)=0$ provided $\mu>\nu$.
\end{lemma}

Finally, we introduce the key notion of this section.

\begin{definition}\label{dhn} Let $X$ be a representation of $Q$. A 
Harder-Narasimhan (HN) 
filtration of $X$ is a filtration $0=X_0\subset X_1\subset\ldots\subset X_s=X$
such that the quotients $X_k/X_{k-1}$ are semistable for $k=1\ldots s$, and
$\mu(X_1/X_0)>\mu(X_2/X_1)>\ldots >\mu(X_s/X_{s-1})$.
\end{definition}

\begin{proposition}\label{ehnf} Any ${\bf k}$-representation $X$ of $Q$ 
posesses a 
unique 
Har\-der-Na\-ra\-sim\-han 
filtration.
\end{proposition}

\proof We proceed by induction on the dimension of $X$. Let $X_1\subset X$ be 
a subrepresentation of maximal slope, and of maximal dimension 
among 
the subrepresentations with this property. It is easy to see that this 
determines $X_1$ uniquely. By induction, $X/X_1$ possesses a unique HN 
filtration, which we can lift to one of $X$ via the projection $X\rightarrow 
X_1$. Uniqueness follows inductively from the fact that the first term of a HN 
filtration already has to be the $X_1$ just constructed. \hb

\begin{definition} An element $d\in{\bf N}I$ is called semistable if there 
exists a se\-mi\-sta\-ble representation of dimension type $d$. A tuple 
$d^*=(d^1,\ldots,d^s)$ is
called a HN type if each $d^k$ is semistable and $\mu(d^1)>\ldots>\mu(d^s)$.
The sum $|d^*|=\sum_{k=1}^sd^k$ is called the weight of $d^*$, and $l(d^*)=s$ 
is called the length of $d^*$.
\end{definition}

\section{Representation varieties and the HN stratification}

We first recall the varieties of quiver representations of fixed dimension 
type. 

\begin{definition} For $d=\sum_{i\in I}d_ii\in{\bf N}I$, define 
$R_d=\bigoplus_{\alpha:i\rightarrow j}{\rm Hom}_{\bf k}({\bf k}^{d_i},{\bf 
k}^{d_j})$
and
$G_d=\prod_{i\in I}{\rm GL}_{d_i}({\bf k})\subset E_d=\bigoplus_{i\in I}{\rm 
End}({\bf k}^{d_i})$. The group $G_d$ acts on $R_d$ via
$$(g_i)_i\cdot(X_\alpha)_\alpha=(g_jX_\alpha g_i^{-1})_{\alpha:i\rightarrow 
j}.$$
\end{definition}

$R_d$ is an affine ${\bf k}$-variety parametrizing the ${\bf 
k}$-representations of $Q$
of dimension type $d$. The $G_d$-orbits ${\cal O}_X$ in $R_d$ correspond
bijectively to the isomorphism classes $[X]$ of representations 
$X\in\bmod_{\bf k}Q$ 
of dimension type $d$. The unique closed orbit of $G_d$ in $R_d$ corresponds 
to the unique semisimple representation of dimension type $d$, so the 
invariant ring ${\bf k}[R_d]^{G_d}$ reduces to the scalars.\\[1ex]
The main use of the notion of (semi-)stability lies in the following result.

\begin{theorem}\cite{Ki} Let ${\bf k}$ be algebraically closed. For each 
$d\in{\bf 
N}I$, the subset $R_d^{ss}\subset R_d$ of semistable representations is an 
open subvariety. It admits a categorial (GIT) quotient ${\cal 
M}^{ss}_d=R_d^{ss}/G_d$, which is a projective variety. The quotient ${\cal 
M}^{ss}_d$ contains a smooth open subvariety ${\cal M}^s_d$, which is a 
geometric quotient by $G_d$ of $R_d^{s}\subset R_d^{ss}$, the subset of stable 
representations.
\end{theorem}

\remark In fact, the notions of (semi-)stability for representations 
translate into the corresponding notions of \cite{Mu} for points in $R_d$, 
with respect to the trival line bundle on $R_d$ with $G_d$-action twisted by 
a suitable character $\chi$ of $G_d$. The 
moduli space 
can be defined as
$${\cal M}^{ss}_d={\bf Proj}(\oplus_{n\in{\bf N}}{\bf k}[R_d]^{G_d,\chi^n}),$$
where 
${\bf k}[R_d]^{G_d,\chi^n}$ consists of semi-invariant polynomial functions 
$f:R_d\rightarrow {\bf k}$ of weight $\chi^n$.\\[2ex]
\examples {}
\begin{enumerate}
\item In example B of section \ref{hnf}, we recover the quotient $({\bf 
P}_{m-1})^n_{\rm stable}/{\rm PGL}_m$ of (\cite{Mu}, 3).
\item The situation of example D of section \ref{hnf} is related to moduli of 
bundles on ${\bf P}^2$, see \cite{D}.
\item As a particular case of example D, the moduli space ${\cal 
M}_{mi+mj}^{ss}$ compactifies the 
affine variety given by ${\bf k}[M_m({\bf k})^{n-1}]^{{\rm GL}_m}$, the ring 
of 
invariants of 
$n-1$-tuples of $m\times m$-matrices under simultaneous conjugation. In fact, 
if the first map $f_1$ is invertible, then its stabilizer $G\simeq{\rm GL}_m$ 
under the ${\rm 
GL}_m\times {\rm GL}_m$-action acts by simultaneous conjugation on 
$f_2,\ldots, f_n$.
\end{enumerate}

Next, we introduce the Harder-Narasimhan stratification. Note that the term 
stratification is used in a weak sense, meaning a finite decomposition of a 
variety into irreducible, locally closed subsets (see \cite{Kir}, Introduction).

\begin{definition} For a HN type $d^*$, we denote by $R_{d^*}^{HN}\subset R_d$ 
the subset of
representations whose HN filtration is of type $d^*$.
$R_{d^*}^{HN}$ is called the HN stratum for the HN type $d^*$. More generally, 
we denote by $R_d^{d^*}$ the subset of
representations possessing some filtration of type $d^*$, i.e.~$\dimv 
X_k/X_{k-1}=d^k$ for $k=1\ldots s$ in the notation of Definition \ref{dhn}.
\end{definition}

\begin{proposition}\label{hnstr} The HN strata for HN types $d^*$ of weight 
$d$ define a stratification of $R_d$ into irreducible, locally closed
subvarieties. The codimension of $R_{d^*}^{HN}$ in $R_d$ is given by 
$-\sum_{1\leq 
k<l\leq s}\langle d^k,d^l\rangle$.
\end{proposition}

\proof Let $F^*\, :\, 0=F^0\subset F^1\subset\ldots\subset F^s$ be a flag of 
type $d^*$ in the $I$-graded vector space ${\bf k}^d=\oplus_{i\in I}{\bf 
k}^{d_i}$, i.e.~$F^k/F^{k-1}\simeq {\bf k}^{d^k}$ 
for $k=1\ldots s$, and denote by $F_i^k$ the $i$-component of $F^k$.
Denote by $Z$ the closed subvariety of $R_d$ of 
representations $X$ which are compatible with $F^*$, i.e.~$X_\alpha(F_i^k)\subset F_j^k$ for $k=1\ldots s$ and for all arrows
$\alpha:i\rightarrow j$ in $Q$. It  is easy 
to see that $Z$ is a trivial vector bundle of rank 
$\sum_{k<l}\sum_{i\rightarrow 
j}d^l_id^k_j$ over $R_{d^1}\times\ldots\times R_{d^k}$, via the projection $p$ 
mapping $X\in Z$ to the sequence of subquotients with respect to $F^*$. Thus, 
the inverse image of $R^{ss}_{d^1}\times\ldots\times R^{ss}_{d^s}$ under $p$ 
is an open subvariety $Z_0$ of $Z$. The action of $G_d$ on $R_d$ induces 
actions of the parabolic subgroup $P_{d^*}$ of $G_d$, consisting of elements
fixing the flag $F^*$, on $Z_0$ and $Z$. The image of the associated fibre 
bundle $G_d\times^{P_{d^*}}Z$ under the action morphism $m$ equals 
$R_d^{d^*}$, 
which is thus a closed subvariety of $R_d$. The image of 
$G_d\times^{P_{d^*}}Z_0$ under $m$ equals $R_{d^*}^{HN}$, and $G_d\times^{P_{d^*}}Z_0$ is the full preimage. By the 
uniqueness 
of the HN filtration, the morphism $m$ is bijective over $R_{d^*}^{HN}$, which 
therefore is an irreducible, open subvariety of $R_d^{d^*}$. The codimension 
is now easily computed as $-\sum_{k<l}\langle d^k,d^l\rangle$, using the 
identity $\langle d,d\rangle=\dim G_d-\dim R_d$ and the above description of 
$Z$. \hb

From this description, we can derive a recursive criterion for the existence 
of semistable
representations, complementing a criterion for the existence of stable 
representations in \cite{Ki}. 

\begin{corollary}\label{ess} A dimension type $d\in{\bf N}I$ is semistable 
if and only
if there exists no HN type $d^*$ of weight $d$ such that $\langle d^k,d^l
\rangle=0$ for all $1\leq k<l\leq s$.
\end{corollary}

\proof A dimension type $d$ is semistable if and only if there exists no 
dense HN stratum, i.e.~a stratum $R_{d^*}^{HN}$ of codimension $0$, which by 
the above codimension formula means $\langle d^k,d^l
\rangle=0$ for all $1\leq k<l\leq s$. \hb

We have seen in the proof of Proposition \ref{hnstr} that the closure of the 
HN stratum $R_{d^*}^{HN}$ equals $R_d^{d^*}$. In general, this is not anymore 
a union of HN strata (see the example below). But at least it is contained in 
such a union, and the strata involved in this union can be controlled using 
polygonal combinatorics. We proceed as in the vector bundle situation, 
following \cite{Sh}.

\begin{definition} Given an arbitrary tuple $d^*=(d^1,\ldots,d^s)$ in ${\bf 
N}I$, we denote by $P(d^*)$ the polygon in ${\bf N}^2$ with
vertices $(\sum_{l=1}^k\dim d^l,\sum_{l=1}^k\Theta(d^l))$ for $k=0\ldots s$. 
For two such tuples
$d^*,e^*$, define $d^*\leq e^*$ (resp.~$d^*<e^*$) if $P(d^*)$ lies on or below 
(resp.~strictly below) $P(e^*)$ (see 
\cite{Sh}). We call $d^*$ convex if the polygon 
$P(d^*)$ is convex.
\end{definition}

Note that, by definition, HN types are always convex.

\begin{proposition} The closure $R_d^{d^*}$ of the HN stratum for the HN type 
$d^*$ is
contained in the union of the HN strata for HN types $e^*\geq
d^*$.
\end{proposition}

\proof Let $X$ be a representation in $R_d^{d^*}$, and  let $e^*$ be the HN 
type of $X$. Since the polygon $P(e^*)$ is convex by definition, it suffices 
to prove that the slope $\mu(\dimv U)$ of an arbitrary subrepresentation
$U\subset X$ lies on or below $P(e^*)$. We proceed by induction over the 
length $t$ of $e^*$; in case $t=1$ there is nothing to prove, since then $X$ 
is semistable. By induction, the slope of $(U+X_1)/X_1$ lies on or below 
$P(e^{*+1})$, where $e^{*+1}$ is the HN type of $X/X_1$. It follows that the 
slope of $U+X_1$ lies on or below $P(e^*)$. Using the obvious exact sequences 
relating $U$, $X_1$ and $U+X_1$, we get by the maximality property of $X_1$:
$$\mu(U+X_1)\leq\mu(X_1)\leq\mu(X_1/(U\cap X_1))\mbox{, thus 
}\mu(U)\leq\mu(U+X_1),$$
which proves the desired property. \hb

\example Let $Q$ be the quiver $i\rightarrow j\rightarrow k$, $d=i+j+k$, 
$\theta=2i^*+3j^*$. Thus $R_d\simeq {\bf k}^2$, where a point $(x,y)$ 
corresponds to 
the representation ${\bf k}\stackrel{x}{\rightarrow}{\bf 
k}\stackrel{y}{\rightarrow}{\bf k}$.
The HN stratum for the HN type $(i,j+k)$ equals ${\bf k}^*\times 0\subset {\bf 
k}^2$, and 
the HN stratum for the HN type $(j,i,k)$ equals $0\times {\bf k}\subset {\bf 
k}^2$. Thus, 
the closure of $R_{(i,j+k)}$ is contained in $R_{(i,j+k)}\cup R_{(j,i,k)}$, 
but is not a union of HN strata.

\section{Quantum groups and the HN system}\label{hns}

In this and the following sections, ${\bf k}$ denotes a finite field, whose cardinality is denoted by $v^2$ for $v\in{\bf
C}$.\\[1ex] To get a relation between the methods of the previous sections and (quantized) 
Kac-Moody 
algebras, we first define the Hall algebra (see \cite{Ri}, \cite{Gr}) in the 
version of Kapranov (see \cite{Ka}).

\begin{definition}\label{hall} For $d\in{\bf N}I$, define ${\cal H}_d={\bf 
C}^{G_d}[R_d]$, 
the ${\bf C}$-vector space of (arbitrary) $G_d$-invariant functions from 
$R_d$ to ${\bf 
C}$.
Define a ${\bf N}I$-graded ${\bf C}$-vector space ${\cal H}={\cal H}(Q)=\bigoplus_{d\in{\bf N}I}{\cal H}_d$
with the multiplication
$$(f*g)(X)=v^{\langle e,d\rangle}\sum_{U\subset X}f(U)g(X/U)$$
on homogeneous elements $f,g$ of degree $|f|=d$ (i.e.~$f\in{\cal H}_d$) and $e$, respectively.
Define a scalar product on ${\cal H}$ by
$$(f,g)=(\# G_d)^{-1}\sum_{X}f(X)g(X)$$
for $f,g\in{\cal H}_d$, where $\#S$ denotes the cardinality of a finite set 
$S$, and $({\cal H}_d,{\cal H}_e)=0$ for $d\not=e$.
\end{definition}

By \cite{Ri}, \cite{Gr}, ${\cal H}$ is an associative, ${\bf N}I$-graded ${\bf 
C}$-algebra.\\[1ex]
An immediate induction yields:

\begin{lemma}\label{multfilt} For $f_k\in{\cal H}_{d^k}$, $k=1\ldots s$, we 
have
$$(f_1*\ldots*f_s)(X)=v^{\langle d^*\rangle}\sum_{0=X_0\subset 
X_1\subset\ldots\subset X_s=X}f_1(X_1/X_0)\ldots f_s(X_s/X_{s-1}),$$
where
$$\langle d^*\rangle=\sum_{k<l}\langle d^l,d^k\rangle.$$
\end{lemma}

Next, we define the so-called composition subalgebra of ${\cal H}$.

\begin{definition} For $d\in{\bf N}I$, define $\chi_d=\chi_{R_d}$, the 
characteristic function of the variety $R_d$. Define ${\cal C}$ as the ${\bf 
C}$-subalgebra of
${\cal H}$ generated by the functions $\chi_i$ for $i\in I$.
\end{definition}

\begin{lemma}\label{ca} We have $\chi_d\in{\cal C}$ for all $d\in{\bf N}I$.
\end{lemma}

\proof We order the set of vertices of $Q$ as $I=\{1,\ldots,n\}$ such that 
$i>j$ if there exists an arrow $i\rightarrow j$, which is possible since $Q$ 
has no oriented cycles. It is easy to see that any representation $X\in R_d$ 
has a unique filtration $0=X_0\subset\ldots\subset X_n=X$ such that $\dimv 
X_{k}/X_{k-1}=d_kk\in{\bf N}I$, so that the product 
$\chi_{d_11}*\ldots*\chi_{d_nn}$ 
equals $\chi_d$ up to some power of $v$. But since $R_{ni}$ for $n\in{\bf 
N}$, $i\in I$ consists of a single point, $\chi_{ni}$ equals a (non-zero) 
scalar multiple 
of $\chi_i^{*n}$. This shows that $\chi_d$ is generated by the $\chi_i$. \hb

The connection to quantum groups is given by the following theorem, due to
J.~A.~Green \cite{Gr}, which generalizes a theorem of C.~M.~Ringel \cite{Ri} 
in the finite type case. 
\begin{theorem}The composition algebra ${\cal C}$ is isomorphic to ${\cal 
U}_v({\frak 
n}^+)$, the quantized enveloping
algebra of the positive part of the symmetric Kac-Moody algebra corresponding 
to the 
quiver $Q$, specialized 
at $v$.
\end{theorem}

Here, the Kac-Moody algebra corresponding to $Q$ is defined by the generalized 
Cartan matrix $(\langle i,j\rangle+\langle j,i\rangle)_{i,j\in I}$ as in 
\cite{Kac}.\\[2ex]
The link between the HN filtration and the Hall algebra is provided by the
following definition.

\begin{definition} For a HN type $d^*=(d^1,\ldots,d^s)$ of weight $d$, define
$$\chi^{HN}_{d^*}=\chi_{R_{d^*}^{HN}}\in{\cal H}_d,$$
the characteristic function of the HN stratum for $d^*$. In particular, define
$\chi^{ss}_{d}=\chi^{HN}_{(d)}$.
\end{definition}

We find the following analogue of the description of the HN stratum in 
Proposition \ref{hnstr}.

\begin{lemma}\label{dechn} For a HN type $d^*=(d^1,\ldots,d^s)$, we have in 
${\cal H}$ the following identity:
$$\chi^{HN}_{d^*}=v^{-\langle 
d^*\rangle}\chi^{ss}_{d^1}*\ldots*\chi^{ss}_{d^s}.$$
\end{lemma}

\proof Any filtration $X_*$ of type $d^*$ with semistable subquotients is 
 a HN filtration by definition, which is unique by Proposition 
\ref{ehnf}. The statement follows from the description Lemma \ref{multfilt} of 
the convolution 
product in ${\cal H}(Q)$. \hb

We easily derive the Harder-Narasimhan recursion, in analogy to \cite{AB} and 
\cite{HN}. 

\begin{proposition}\label{rec} For all $d\in{\bf N}I$, we have
$$\chi^{ss}_d=\chi_d-\sum_{d^*}v^{-\langle 
d^*\rangle}\chi^{ss}_{d^1}*\ldots*\chi^{ss}_{d^s},$$
where the sum runs over all HN types $d^*\not=(d)$ of weight $d$.
We have $\chi^{ss}_d=\chi_d$ if and only if $\Theta$ is constant
on $\supp d=\{i\in I\, :\, d_i\not=0\}$, the support of $d$ in $I$.
\end{proposition}

\proof In case $\Theta$ is constant on $\supp d$, any representation $X\in 
R_d$ is obviously semistable, thus $\chi^{ss}_d=\chi_d$. Otherwise, we have a 
disjoint union $R_d=\bigcup_{d^*}R_{d^*}^{HN}$ (the union being over all HN 
types 
of weight $d$), which translates into the identity 
$\chi_d=\sum_{d^*}\chi_{d^*}^{HN}$. Using Lemma \ref{dechn}, the recursion is 
proved. \hb

As an immediate consequence, we get:

\begin{theorem}\label{hnsystem} For all HN types $d^*$, the element
$\chi^{HN}_{d^*}$ already belongs to the composition algebra ${\cal C}$.
The scalar product of two such elements is given by
$$(\chi^{HN}_{d^*},\chi^{HN}_{e^*})=\delta_{d^*,e^*}(\# 
G_d)^{-1}\# R_{d^*}^{HN}.$$
If $\mu(d^s)>\mu(e^1)$, then $\chi^{HN}_{d^*}*\chi^{HN}_{e^*}=
\chi^{HN}_{(d^1\ldots d^se^1\ldots e^t)}$.
\end{theorem}

\proof The first statement follows from Proposition \ref{rec} by induction on 
the dimension type, using Lemma \ref{ca}. The second statement just 
reformulates the disjointness of different HN strata, using the definition of 
the scalar product on ${\cal H}$. For the third part, we just have to note 
that, by assumption, the concatenation $(d^1,\ldots,d^s,e^1,\ldots,e^t)$ is 
again a HN type. \hb

\begin{definition} The set of elements $\chi^{HN}_{d^*}$ for various HN types
$d^*$ is called the Harder-Narasimhan system in ${\cal C}\simeq{\cal 
U}_v({\frak n}^+)$.
\end{definition}

We thus have found a natural orthogonal system of elements in ${\cal 
U}_v({\frak n}^+)$, which is neccessarily linearly independent by 
orthogonality and the fact that $\chi_{d^*}^{HN}$ is non-zero by definition of 
a HN type. This system is recursively computable by combining 
Lemma \ref{dechn}, Proposition \ref{rec} and Corollary \ref{ess}. Furthermore, 
it is partially multiplicative.
The HN system can be viewed as a replacement for PBW type bases (which 
have 'good' properties only for finite type cases) in the Kac-Moody 
situation.\\[2ex]
We end this section with a compact reformulation of the HN recursion \ref{rec}.

\begin{definition} Consider the skew Laurent polynomial ring $\widehat{\cal 
H}={\cal H}[[T_i, :\, i\in I]]$, where $\chi_i*T_j=v^{\langle i,j\rangle}
T_j*\chi_i$ for all $i,j\in I$.
Define generating functions
$$X(T)=\sum_{d\in{\bf N}I}\chi_d*T^d,\;\; X^{ss}_\mu(T)=\sum_{{d\in{\bf 
N}I}\atop{\mu(d)=\mu}}\chi^{ss}_d*T^d,$$
where $T^d=\prod_{i\in I}T_i^{d_i}$ for $d\in{\bf N}I$.
\end{definition}

Using these definitions, the recursion \ref{rec} allows the following
notation:

\begin{proposition} In $\widehat{\cal H}$, the generation function $X(T)$ 
equals the 'descending product' $\prod^{\leftarrow}_{\mu\in{\bf 
Q}}X^{ss}_\mu(T):=\sum_{\mu_1>\ldots>\mu_s}X_{\mu_1}^{ss}(T)*\ldots *
X_{\mu_s}^{ss}(T)$.
\end{proposition}

\proof First note that for each $f\in{\cal C}\cap{\cal H}_d$, we have 
$f*T^e=v^{\langle d,e\rangle}T^e*f$ by definition. We calculate using Lemma 
\ref{dechn}:
\begin{eqnarray*}
\prod^{\leftarrow}_{\mu\in{\bf 
Q}}X^{ss}_\mu(T)&=&\sum_{\mu_1>\ldots>\mu_s}X_{\mu_1}^{ss}(T)*\ldots *
X_{\mu_s}^{ss}(T)\\
&=&\sum_{d^*}\chi_{d^1}^{ss}*T^{d^1}*\ldots*\chi_{d^s}^{ss}*T^{d^s}\\
&=&\sum_{d^*}v^{-\sum_{k<l}\langle 
d^l,d^k\rangle}\chi_{d^1}^{ss}*\ldots*\chi_{d^s}^{ss}*T^{d^1}*\ldots *T^{d^s}\\
&=&\sum_{d^*}\chi_{d^*}^{HN}*T^{|d^*|}=\sum_{d\in{\bf NI}}\chi_d*T^d=X(T),
\end{eqnarray*}
where the sums run over all HN types $d^*$.\hb

\section{Resolving the recursion}

Motivated by similar results in \cite{Za},\cite{LR} (but with different 
methods), we will now derive a resolution of the recursion \ref{rec}. This
gives an explicit formula for the
elements $\chi^{ss}_d$:

\begin{theorem}\label{ex} For all $d\in{\bf N}I$, we have:
$$\chi^{ss}_d=\sum_{d^*}(-1)^{s-1}v^{-\langle d^*\rangle}\chi_{d^1}*\ldots
*\chi_{d^s},$$
where the sum runs over all tuples of non-zero dimension types $d^*=(d^1\ldots 
d^s)$ of weight $d$
such that $d^*=(d)$ or $d^*>(d)$, i.e.~$\mu(\sum_{l=1}^kd^l)>\mu(d)$ for 
$k=1\ldots 
s-1$.
\end{theorem}

\proof The proof proceeds along the following lines:\\[1ex]
First, we reduce the statement 
by an explicit calculation to a purely combinatorial formula. This will be 
interpreted in terms of the Euler characteristic of a simplicial complex 
(encoding convex coarsenings of a polygon). By explicit combinatorics, we show 
that this simplicial 
complex is in fact always contractible, proving the Theorem.\\[1ex]
We start with some definitions. 

\begin{definition} Let $d^*=(d^1,\ldots,d^s)$ be a tuple of dimension types.
\begin{enumerate}
\item For a subset $I=\{s_1<\ldots<s_k\}\subset\{1,\ldots,s-1\}$,
define the $I$-coarsening of $d^*$ as
$$c_I^*(d^*)=(d^1+\ldots+d^{s_1},d^{s_1+1}+\ldots+d^{s_2},\ldots,d^{s_k+1}+\ldots+d^s).$$
\item The subset $I$ is called $d^*$-admissible if the following holds:
\begin{enumerate}
\item $c_I^*(d^*)$ is convex,
\item For all $i=0\ldots k$, we have $(d^{s_i+1},\ldots,d^{s_{i+1}})\geq 
(d^{s_i+1}+\ldots+d^{s_{i+1}})=(c_I^i(d^*))$.
\end{enumerate}
\item Define ${\cal A}(d^*)$ as the set of all $d^*$-admissible subsets of 
$\{1,\ldots,s-1\}$.
\end{enumerate}
\end{definition}

In polygonal language, a coarsening is thus admissible if it is convex and 
lies on or below the polygon $P(d^*)$.

\begin{lemma} For all sequences $d^*$, the set ${\cal A}(d^*)$ is a simplicial 
complex.
\end{lemma}

\proof If $I$ is $d^*$-admissible, and $J\subset I$ is a subset, then 
$c_J^*(d^*)=c_J^*(c_I^*(d^*))$. Thus, $c_J^*(d^*)$ is a coarsening of the 
convex tuple $c_I^*(d^*)$. Therefore, it obviously has to be convex again, and 
it has to lie on or below $c_I^*(d^*)$, which lies on or below $P(d^*)$ by 
assumption. We 
conclude that $J$ is also $d^*$-admissible. \hb

Using these definitions, we can reduce the formula of Theorem \ref{ex} to a 
purely combinatorial problem. We adopt Proposition \ref{rec}. If $\Theta$ is 
constant on $\supp d$, then any dimension type $e\leq d$ has the same slope as 
$d$, which means that the sum in Theorem \ref{ex} just runs over the single 
tuple $(d)$, and the formula is trivial.\\[1ex]
Otherwise, we consider the HN recursion and replace any term $\chi_{d^i}^{ss}$ 
on its right hand side by the claimed formula:
\begin{eqnarray*}
\chi_d&=&\sum_{d^*}v^{-\langle 
d^*\rangle}\chi_{d^1}^{ss}*\ldots*\chi_{d^s}^{ss}\\
&=&\sum_{d^*}\sum_{d^{1,*},\ldots,d^{s,*}}v^{-\langle 
d^*\rangle}(-1)^{\sum_{i=1}^s(t_i-1)}v^{-\sum_{i=1}^s\langle 
d^{i,*}\rangle}\times\\
&&\times 
\chi_{d^{1,1}}*\ldots\chi_{d^{1,t_1}}*\ldots*\chi_{d^{s,1}}*\ldots*\chi_{d^{s,t_s}}.
\end{eqnarray*}
In this equation, the outer sum runs over all HN types $d^*$ of weight $d$, 
and the inner sums run over all sequences $(d^{1,*},\ldots,d^{s,*})$ of 
tuples $d^{i,*}$ of weight $d^i$ and length $t_i$ such that 
$d^{i,*}=(d^i)$ or $d^{i,*}>(d^i)$.\\[1ex]
We now want to exchange the order of summation. So denote by $e^*$ the 
concatenation $e^*=(d^{1,1},\ldots,d^{s,t_s})$ of all the tuples $d^{i,*}$. 
Obviously, the resulting tuples $e^*$ run over all tuples of weight $d$ such 
that $e^*=(d)$ or $e^*>(d)$. By the above definitions, the HN type $d^*$ of 
the outer 
sum is an admissible coarsening of $e^*$. The sign and the $v$-exponent in the 
above sum are easily computed as $(-1)^{l(e^*)-l(d^*)}v^{-\langle 
e^*\rangle}$. So the above sum can be rewritten as:
\begin{eqnarray*}
\chi_d&=&\sum_{e^*}(-1)^{l(e^*)-1}v^{-\langle 
e^*\rangle}\sum_{d^*}(-1)^{l(d^*)-1}\chi_{e^1}*\ldots*\chi_{e^{l(e)}},
\end{eqnarray*}
where the outer sum runs over all tuples $e^*$ of weight $d$ such that 
$e^*=(d)$ or $e^*\geq d$, and the inner sum runs over all admissible 
coarsenings of $e^*$. 
Theorem \ref{ex} thus reduces to the following statement:

\begin{lemma}\label{ossa} For each tuple $d^*$ of weight $d$ such that 
$d^*=(d)$ or 
$d^*>(d)$, we 
have
$$\sum_{I\in{\cal A}(d^*)}(-1)^{\# 
I}=\left\{\begin{array}{ccc}0&,&d^*\not=(d),\\
1&,&d^*=(d)\end{array}\right.$$
\end{lemma}

\proof We proceed by induction on $s=l(d^*)$. In case $s=1$, we have 
$d^*=(d)$, and there is nothing to prove. Otherwise, consider the slopes of 
the first two entries $d^1$, $d^2$.\\[1ex]
If $\mu(d^1)<\mu(d^2)$, we define $I_0\subset\{1,\ldots,s-1\}$ as the subset 
of all $k$ such that $\mu(d^1)\geq\mu(d^1+\ldots+d^k)$. It is then easy to see 
from the definitions that ${\cal A}(d^*)={\cal A}(c_{I_0}^*(d^*))$. Since 
$2\not\in I_0$, the Lemma holds for the $I_0$-coarsening by induction.\\[1ex]
If $\mu(d^1)\geq\mu(d^2)$, we define $d'^*=(d^2,\ldots,d^s)$ and 
$d''^*=(d^1+d^2,d^3,\ldots,d^s)$. Again, it is then easy to see from the 
definitions that we have a disjoint union
$${\cal A}(d^*)={\cal A}(d''^*)\cup(\{1\}\cup({\cal A}(d'^*)+1)),$$
where $\{1\}\cup({\cal A}(d'^*)+1)$ consists of the sets 
$\{1,s_1+1,\ldots,s_k+1\}$ for 
$\{s_1,\ldots,s_k\}\in{\cal A}(d'^*)$. This leads to the calculation
$$\sum_{I\in{\cal A}(d^*)}(-1)^{\# I}=\sum_{I\in{\cal A}(d''^*)}(-1)^{\# 
I}+\sum_{I\in{\cal A}(d'^*)}(-1)^{\# I+1},$$
which thus equals zero by induction. \hb

Combining the above calculation and the Lemma, we see that 
Theorem \ref{ex} is proved.\hb\\[3ex]
{\remarks $\;$
\begin{enumerate}
\item Implicitely in the proof of the above Lemma, we have proved the 
contractibility of the simplicial complex ${\cal A}(d^*)$.
\item The statement is easily generalized to give an explicit formula for 
arbitrary
elements of the HN system, using Lemma \ref{dechn}.
\item The summation can be viewed as running over all lattice paths with 
arbitrary step length in the lattice
of all elements $e\in{\bf N}I$ such that $e\leq d$ and $\mu(e)>\mu(d)$.
\end{enumerate}}

The last remark gives a key to a compact reformulation of the above formula. 
We adopt a version of the transfer matrix method (see \cite{St}).

\begin{corollary}\label{tma} Let ${\cal T}_d$ be the quadratic matrix with 
rows and
columns indexed by ${\cal I}(d)=\{e\in{\bf N}I\, :\, e\leq d,\, 
\mu(e)>\mu(d)\}\cup\{0,d\}$, and with entries in ${\cal C}$
given by $v^{\langle e-f,e\rangle}\chi_{f-e}$ if $e\leq f$, and zero otherwise.
Then the $(0,d)$-entry of the inverse matrix ${\cal T}_d^{-1}$ equals 
$-\chi^{ss}_d$.
\end{corollary}

\proof Since the matrix ${\cal T}$ is upper unitriangular with respect to the 
ordering $\leq$ on ${\cal I}(d)$, an entry of the inverse matrix can be 
computed as
\begin{eqnarray*}
-({\cal T}^{-1})_{0,d}&=&\sum_{0=e^0<e^1<\ldots<e^s=d}(-1)^{s-1}{\cal 
T}_{e^0,e^1}\ldots{\cal T}_{e^{s-1},e^s}\\
&=&\sum_{0=e^0<e^1<\ldots<e^s=d}(-1)^{s-1}v^{\langle 
e^0-e^1,e^0\rangle+\ldots+\langle e^{s-1}-e^s,e^{s-1}\rangle}\times\\
&&\times\chi_{e^1-e^0}*\ldots*\chi_{e^s-e^{s-1}}\\
&=&\sum_{d^1,\ldots,d^s}(-1)^{s-1}v^{-\sum_{k<l}\langle 
d^l,d^k\rangle}\chi_{d^1}*\ldots*\chi_{d^s}
\end{eqnarray*}
by substituting $d^k=e^k-e^{k-1}$. Noting that the defining properties of 
${\cal I}(d)$ translate into $\mu(d^k)>\mu(d)$ for all $k=1\ldots s-1$, we 
arrive at the formula of Theorem \ref{ex}. \hb

This reformulation of the formula will be used in the next section to derive a 
fast algorithm for the computation of Betti numbers (see Corollary \ref{algo}).

\section{Cohomology of quiver moduli}\label{cqm}

We still keep the assumption that ${\bf k}$ is a finite field with $v^2$ 
elements. 
As the main application of the HN system and its explicit formula Theorem 
\ref{ex}, we derive a formula for the Poincare polynomial of ordinary 
cohomology of the 
complex moduli spaces ${\cal M}_d^{ss}({\bf C})$ 'in the coprime case' (see 
below for the 
definition).\\[1ex]
The strategy is similar to the approach of \cite{HN}: First, we count numbers 
of rational points of varieties. Then, we use results of \cite{Mu} to prove 
compatibilities between these numbers. Finally, we relate these numbers to 
Betti numbers of complex varieties with the aid of the Weil 
conjectures.\\[1ex]
The link between the HN system and numbers of rational points is provided by a 
twisted character of the Hall algebra ${\cal H}$.

\begin{lemma} The map $ev:{\cal H}\rightarrow {\bf C}$, defined by 
$ev(f)=(\# G_d)^{-1}\sum_{X}f(X)$ on ${\cal H}_d$, fulfills
$$ev(f*g)=v^{-\langle e,d\rangle}ev(f)ev(g)$$
for $|f|=d$, $|g|=e$.
\end{lemma}

\proof Without loss of generality, we can assume $f$ (resp.~$g$) to be the 
characteristic function of an orbit ${\cal O}_M$ (resp.~${\cal O}_N$). By the 
definitions, we have $ev(\chi_{{\cal O}_M})=\# {\rm Aut}(M)^{-1}$, and 
$\chi_{{\cal O}_M}*\chi_{{\cal O}_N}=v^{\langle 
e,d\rangle}\sum_{[X]}F_{N,M}^X\chi_{{\cal O}_X}$, where $F_{N,M}^X$ 
denotes the number of subrepresentations of $X$ which are isomorphic to $M$, 
with quotient isomorphic to $N$. By a formula of C. Riedtmann \cite{Rie}, we 
have
$$F_{N,M}^X=v^{-2\dim{\rm Hom}(N,M)}\frac{\#{\rm Aut}(X)}{\#{\rm 
Aut}(M)\cdot\# {\rm Aut}(N)}\#{\rm Ext}^1(N,M)_X,$$
where ${\rm Ext}^1(N,M)_X$ denotes the set of extension classes with middle 
term isomorphic to $X$. Using this formula, the Lemma follows by an easy 
calculation. \hb

Applying the evaluation map $ev$ to both sides of the formula Theorem 
\ref{ex}, using its definition and the previous lemma, we get immediately:

\begin{corollary}\label{crp} For all $d\in{\bf N}I$, we have:
$$\frac{\# R_d^{ss}}{\# G_d}=\sum_{d^*}(-1)^{s-1}v^{-2\langle 
d^*\rangle}\prod_{k=1}^s\frac{\# R_{d^k}}{\# G_{d^k}},$$
where the sum runs over all tuples of non-zero dimension types $d^*=(d^1\ldots 
d^s)$ of weight $d$
such that $\mu(\sum_{l=1}^kd^l)>\mu(d)$ for $k=1\ldots 
s-1$.
\end{corollary}

We thus have to relate the number of points in $R_d^{ss}$ to the number of 
points in the geometric quotient ${\cal M}_d^{ss}$. This is directly possible 
only in 
the following case.

\begin{definition} A dimension type $d\in{\bf N}I$ is called coprime if 
the numbers $\Theta(d),\dim d\in{\bf Z}$ are coprime.
\end{definition}

(This property seems to be very restrictive at first sight; see however 
section \ref{examples} for enough interesting examples.)

\begin{lemma}\label{l1} For coprime $d$, we have $R_d^{ss}=R_d^{s}$, and 
${\rm End}_{{\bf k}Q}(X)\simeq {\bf k}$ for all $X\in R_d^{ss}$.
\end{lemma}

\proof If $X\in R_d^{ss}$ is not stable, there exists a proper 
subrepresentation $U\subset X$ such that ${\Theta(U)}/{\dim 
U}={\Theta(X)}/{\dim X}$, contradicting coprimality. Consider now the 
scalar 
extension $\overline{X}=\overline{{\bf k}}\otimes_{\bf k}X$ to the algebraic 
closure of ${\bf k}$; we claim that it 
is still 
semistable. In fact, its HN filtration is unique, hence stable under 
Frobenius, and thus it descends to the HN filtration of $X$, which is trivial 
by semistability of $X$. We conclude that $\overline{X}$ is semistable, having 
trivial HN filtration. By the first part of the Lemma, $\overline{X}$ is 
already stable, thus its endomorphism ring reduces to scalars by a Schur's 
Lemma type argument. Thus the same holds for $X$. \hb

\begin{lemma}\label{l2} For coprime $d$, the action of $PG_d=G_d/{\bf k}^*$ on 
$R_d^{ss}$ is free in the sense of Mumford (see (\cite{Mu}, 0.8. iv)).
\end{lemma}

\proof Since $PG_d$ acts set-theoretically free on $R_d^{ss}$ by Lemma 
\ref{l1}, the natural map $\Psi:PG_d\times R_d^{ss}\rightarrow R_d^{ss}\times 
R_d^{ss}$ is injective; we have to prove that it is a closed immersion. 
Consider the map $\Phi:R_d^{ss}\times R_d^{ss}\rightarrow {\rm 
Hom}_{\bf k}(E_d,R_d)$ given by
$$\Phi(X,Y)(\phi_i)_{i\in 
I}=(\phi_jX_\alpha-Y_\alpha\phi_i)_{\alpha:i\rightarrow j}.$$
Since the kernel of $\Phi(X,Y)$ can be identified with the space of 
${\bf k}Q$-ho\-mo\-mor\-phisms from $X$ to $Y$, the image $Z$ of $\Psi$ is 
precisely 
the set of pairs $(X,Y)$ such that $\Phi(X,Y)$ has non-trivial kernel (i.e.~
where $\Phi(X,Y)$ has a fixed rank $r$) by Lemma \ref{l1}. On the open subset 
of $Z$ where a fixed $r\times r$ minor of $\Phi(X,Y)$ is non-vanishing, we can 
thus recover from the pair $(X,Y)$ a matrix $\not=0$ in $E_d$ intertwining $X$ 
and 
$Y$, i.e.~we can 
recover algebraically the unique element of $PG_d$ mapping $X$ to $Y$. Thus, 
we have constructed 
locally an inverse morphism $\Psi^{-1}:Z\rightarrow PG_d\times R_d^{ss}$. \hb

\begin{proposition}\label{nrp} For coprime $d$, we have the following formula 
for the number of ${\bf k}$-rational points:
$$\# {\cal M}_d^{ss}=\frac{\# R_d^{ss}}{\# PG_d}.$$
\end{proposition}

\proof From the proof of Lemma \ref{l1} above we see that semistability is stable under base change. By coprimality, the
same holds for stability. Together with Proposition 1.14. of \cite{Mu}, this implies that
$R_d^{ss}$, the variety of semistable representations over ${\bf k}$, coincides with the semistable locus of the ${\bf k}$-variety
$R_d$. Thus, there exists a uniform geometric quotient 
$\pi:R_d^{ss}\rightarrow {\cal M}_d^{ss}$ (in the sense of (\cite{Mu}, 0.7.)) 
of $R_d^{ss}$ by $PG_d$. Using Lemma \ref{l2}, we can apply (\cite{Mu}, 0.9.) 
to conclude that $\pi$ turns $R_d^{ss}$ into a principal $PG_d$-bundle over 
${\cal M}_d^{ss}$, in the sense that $PG_d\times R_d^{ss}\simeq R_d^{ss}\times_{{\cal M}_d^{ss}}R_d^{ss}$. We conclude that each
fibre of $\pi$ is isomorphic to $PG_d$, and the formula for the number of 
${\bf k}$-rational points follows. \hb

Using the Weil conjectures \cite{Del}, we can conclude:

\begin{theorem}\label{betti} Assume that $d$ is coprime, and let ${\cal 
M}_d^{ss}({\bf C})$ 
be the moduli space of semistable representations of $Q$ over the field 
${\bf k}={\bf C}$. Then the Poincare polynomial of the cohomology with complex 
coefficients of ${\cal M}_d^{ss}({\bf C})$ 
is given by
$$\sum_{i\in{\bf Z}}\dim_{\bf C}H^i({\cal M}_d^{ss}({\bf 
C}))v^i=(v^2-1)\sum_{d^*}(-1)^{s-1}v^{-2\langle 
d^*\rangle}\prod_{k=1}^s\frac{\# R_{d^k}}{\# G_{d^k}},$$
where the sum runs over all tuples of non-zero dimension types $d^*=(d^1\ldots 
d^s)$ of weight $d$
such that $\mu(\sum_{l=1}^kd^l)>\mu(d)$ for $k=1\ldots 
s-1$.
\end{theorem}

\proof By Corollary \ref{crp}, the right hand side equals ${\# 
R_d^{ss}}/{\# G_d}$, which by Proposition \ref{nrp} equals 
$\frac{1}{v^2-1}\#{\cal M}_{d}^{ss}$. A standard argument (see (\cite{Kir}, 
15.), (\cite{Goe}, 1.2.)), together with the generic compatibility of formation of invariants and base change (see (\cite{CBV},
Lemma B.4.)) shows that the number of rational points 
(viewed as a function of $v$) is precisely the Poincare polynomial of
the cohomology of ${\cal M}^{ss}_d({\bf C})$. \hb 

Using the obvious formulas for $\# R_d$ and $\# G_d$, the above formula can be 
simplified and made more explicit. For $N\in{\bf N}$, denote 
$[N]=\frac{v^{2N}-1}{v^2-1}$ and $[N]!=[1][2]\ldots[N]$.

\begin{corollary}\label{exex} For coprime $d=\sum_{i\in I}d_ii$, we have
\begin{eqnarray*}\sum_{i\in{\bf Z}}\dim_{\bf C}H^i({\cal M}_d^{ss}({\bf 
C}))v^i&=&(v^2-1)^{1-\sum_id_i}v^{-\sum_id_i(d_i-1)}\times\\
&&\times\sum_{d^*}(-1)^{s-1}v^{2\sum_{k\leq l}\sum_{i\rightarrow 
j}d^l_id^k_j}\prod_{k=1}^s\prod_i([d^k_i]!)^{-1},
\end{eqnarray*}
where the sum runs over all tuples of non-zero dimension types $d^*=(d^1\ldots 
d^s)$ of weight $d$
such that $\mu(\sum_{l=1}^kd^l)>\mu(d)$ for $k=1\ldots 
s-1$.
\end{corollary}

We can also apply the evaluation map $ev$ to the transfer matrix analogue 
\ref{tma} of the resolved recursion. We easily get:

\begin{corollary}\label{algo} Let $T_d$ be the quadratic matrix with rows and
columns indexed by ${\cal I}(d)$, and with entries in ${\bf C}(v)$
given by $v^{2\langle e-f,e\rangle}\frac{\# R_{f-e}}{\# G_{f-e}}$ if $e\leq 
f$, 
and zero otherwise.
Then for coprime $d$, the $(0,d)$-entry of the inverse matrix $T_d^{-1}$ 
equals 
$(1-v^2)^{-1}\sum_{i\in{\bf Z}}\dim_{\bf C}H^i({\cal M}_d^{ss},{\bf C})v^i$.
\end{corollary}

This last corollary gives a simple and fast algorithm for computing the 
Poincare polynomials, since we just have to solve a linear equation defined by 
the upper unitriangular matrix $T_d$. This is clearly a problem of polynomial 
order. More precisely, it is of quadratic order in the size of the set ${\cal 
I}(d)$, which is approximately quadratic in the product $\prod_i(d_i+1)$. In 
contrast to this, both the HN recursion \ref{rec} and the explicit formula 
\ref{ex} are of exponential order compared to the size of the entries of $d$, 
since the summations run over certain classes of lattice paths.

\section{Applications and examples}\label{examples}

As a first (potential) application of the HN system, we consider the case 
where $Q$ is of Dynkin type, i.e.~the underlying unoriented graph is a 
disjoint 
union of Dynkin diagrams of type $A,D,E$. Generalizing example A of section 
\ref{hnf}, we have the following:

\begin{conjecture} If $Q$ is of Dynkin type, there exists a weight $\Theta$ 
such that the stable representations are precisely the indecomposables.
\end{conjecture}

By Gabriel's theorem (see \cite{ARS}), the dimension types of the 
indecomposables are precisely 
the positive roots for the corresponding root system, and the conjecture 
can be reduced to a purely combinatorial problem.\\[1ex]
Provided the conjecture holds, the HN strata for the corresponding slope 
function $\mu$ are precisely the (finitely many) $G_d$-orbits in $R_d$. The HN 
system is thus a basis for ${\cal U}_v({\frak n}^+)$. By orthogonality, it has 
to coincide with a PBW basis (in the sense of \cite{Lu}) up to scalars. The 
formula \ref{ex} thus gives an explicit description of a PBW basis, and in 
particular of root elements in ${\cal U}_v({\frak n}^+)$.\\[1ex]
Note however that the HN system can never be a basis for infinite
types due to the non-trivial root multiplicities of the corresponding Kac-Moody algebra.\\[2ex]
Another impact of the developed methods is on the structure of the Hall 
algebra ${\cal H}(Q)$ itself, which is much larger than ${\cal C}\simeq{\cal 
U}_v({\frak n}^+)$ if $Q$ is not of Dynkin type. In fact, ${\cal H}$ is a 
specialization of the quantized 
enveloping algebra of a Borcherds algebra by \cite{SV}; the structure of this 
Borcherds algebra remains unknown. Using the concepts of sections \ref{hnf}, 
\ref{hns}, we can define $\mu$-local Hall algebras by ${\cal 
H}_\mu=\bigoplus_{d:\,\mu(d)=\mu}{\bf C}^{G_d}[
R^{ss}_d]$. This is a subalgebra of ${\cal H}$ by Lemma \ref{fl} and 
Definition \ref{hall}. Using 
\cite{SV}, it should again be possible to relate ${\cal H}_\mu$ to a Borcherds 
algebra, whose structure should be intimately related to the geometry of the 
moduli spaces ${\cal M}_d^{ss}$.\\[2ex]
After these conjectural applications, we turn to the examples of section 
\ref{hnf} and make the formula \ref{exex} explicit in some cases.\\[2ex]
In the case of example B, assume that $m$ and $n$ are coprime. Then we can 
apply formula \ref{exex} to this particular case. After some elementary 
reformulations, we get the 
Poincare polynomial of cohomology of the quotient $({\bf P}^{m-1})^n_{\rm 
stable}/{\rm PGL}_m$ as (note the multinomial coefficient):
$$(v^2-1)^{1-m-n}v^{-m(m-1)}\sum_{m_*,n_*}(-1)^{s-1}{n\choose {n_1\ldots 
n_s}}v^{2\sum_{k\leq l}m_kn_l}\prod_k([m_k]!)^{-1},$$
where the sum runs over all tuples $m_*=(m_1\ldots m_s)$, $n_*=(n_1\ldots 
n_s)$ such that $\sum_km_k=m$, $\sum_kn_k=n$, $(m_k,n_k)\not=(0,0)$ for all 
$k$, and $(m_1+\ldots+m_k)/m<(n_1+\ldots+n_k)/n$ for all $k=1\ldots s-1$. This 
formula generalizes the formulas of (\cite{Kir}, 16)\\[2ex]
Similarly, we can deal with example D. Considering the dimension 
vector $ai+bj$ for the quiver $Q=i\stackrel{(n)}{\rightarrow}j$ such that $a$ 
and $b$ are coprime, the Poincare polynomial $P_{a,b}^n(v)$ of the quotient 
$W_{a,b}^n={\rm Hom}({\bf 
C}^a,{\bf C}^b)^n_{\rm stable}/{\rm GL}_a\times {\rm GL}_b$ is given by:
$$(v^2-1)^{1-a-b}v^{-a(a-1)-b(b-1)}\sum_{a_*,b_*}(-1)^{s-1}v^{2n\sum_{k\leq 
l}a_lb_k}\prod_k([a_k]![b_k]!)^{-1},$$
where the sum runs over all tuples $a_*=(a_1\ldots a_s)$, $b_*=(b_1\ldots 
b_s)$ such that $\sum_ka_k=a$, $\sum_kb_k=b$, $(a_k,b_k)\not=(0,0)$ for all 
$k$, and $(a_1+\ldots+a_k)/a>(b_1+\ldots+b_k)/b$ for all $k=1\ldots 
s-1$.\\[1ex]
It is possible to make this formula more tractable; the neccessary 
calculations are elementary, but quite tedious, so they will be omitted here. 
The idea is to 
separate the zero and non-zero entries among the $a_k$, and to use some 
standard identities for the $v$-binomial coefficients $\left[{M\atop 
N}\right]=\frac{[M+N]!}{[M]![N]!}$ (see e.g. \cite{Lu}). The final result is
$$P_{a,b}^n(v)=(v^2-1)^{1-a}v^{-a(a-1)}\sum_{a_*,b_*}(-1)^{s-1}v^{2\sum_{k<l}(na_l-b_l)b_k}
\prod_{k=1}^s([a_k]!)^{-1}\left[{{na_k}\atop{b_k}}\right],$$
where the sum runs over all tuples $a_*=(a_1\ldots a_s)$, $b_*=(b_1\ldots 
b_s)$ such that $\sum_ka_k=a$, $\sum_kb_k=b$, $a_k\not=0$ for all 
$k$, and $(a_1+\ldots+a_k)/a>(b_1+\ldots+b_k)/b$ for all $k=1\ldots 
s-1$.\\[2ex]
Since the summation runs only over non-zero $a_k$, this formula has the 
advantage of being easily computable for small $a$:\\[1ex]
For $a=1$, we just get $\left[{n\atop b}\right]$, the Poincare polynomial of 
the cohomology of the Grassmanian ${\rm Gr}_b^n\simeq W_{1,b}^n$. For $a=2$, 
we get
$$P_{2,b}^n(v)=(v^2-1)^{-1}v^{-2}(\frac{1}{v^2+1}\left[{2n\atop 
b}\right]-\sum_{k=0}^{(b-1)/2}v^{2(n-b+k)k}\left[{n\atop 
k}\right]\left[{n\atop{b-k}}\right]),$$
generalizing results of \cite{D}. This leads to a formula for the Euler 
characteristic (which can not be read off directly from the general formulas):
$$\chi(W_{2,b}^n)=\frac{bn-1}{4}{2n\choose b}-n\sum_{k=0}^{(b-1)/2}k{n\choose 
k}{n\choose b-k},$$
generalizing results of \cite{ES}. Similar results can be obtained for example 
B.\\[2ex]
Finally, let us remark that the algorithm \ref{algo} opens the possibility for 
computer experiments in many non-trivial cases. These experiments suggest 
several formulas for generating functions and asymptotical behaviours of 
Poincare polynomials and Euler characteristics. One may hope that such 
experiments lead to further insights into the geometry of quiver moduli.

\end{document}